\magnification=1200
\baselineskip 15pt
\font\title=cmssdc10 at 15pt

\def\frame #1{\vbox{\hrule height.1pt
\hbox{\vrule width.1pt\kern 10pt
\vbox{\kern 10pt
\vbox{\hsize 3.5in\noindent #1}
\kern 10pt}
\kern 10pt\vrule width.1pt}
\hrule height0pt depth.1pt}}
\def\bframe #1{\vbox{\hrule height.1pt
\hbox{\vrule width.1pt\kern 10pt
\vbox{\kern 10pt
\vbox{\hsize 5in\noindent #1}
\kern 10pt}
\kern 10pt\vrule width.1pt}
\hrule height0pt depth.1pt}}

\footnote{}{Key words and phrases: }
\bigskip
\centerline{\title The age incidence of any cancer can be explained by a one-mutation model }
\bigskip
Rinaldo B. Schinazi

Department of Mathematics

University of Colorado

Colorado Springs, CO 80933-7150

U.S.A.
\bigskip
{\bf Abstract.} We propose a one mutation model for cancer with a mutation rate that increases with time. Under rather general hypotheses the number of mutations is necessarily a (non homogeneous) Poisson process with the prescribed mutation rate. We show that the cumulative probability of cancer up to time $t$ is, up to a multiplicative constant, an antiderivative of the mutation rate. 
\bigskip
{\bf Introduction} Mathematical models explaining age incidence of cancers go back to at least Armitage and Doll (1954). The statistical data show that age incidence follows power laws for several different cancers (with different powers for different cancers). This lead to the hypothesis that a given cancer arises after several successive mutations and that there is a simple relation between the number of necessary mutations and the power law of the incidence. This hypothesis has yet to be proved biologically and several authors  have argued that cancers could be modelled by a two stage process, see Armitage (1985) and Moolgavkar and Knudson (1981). Recently, Michor et al. (2006) proposed a one mutation model to explain the age incidence of chronic myeloid leukemia which increases as a third power of age. In this paper, we show that in fact any power can be explained by a one mutation model. 

In this note, we propose a one mutation model with a mutation rate that increases with time. The number of mutations is necessarily a (non homogeneous) Poisson process with the prescribed mutation rate. We show that the cumulative probability of cancer up to time $t$ is, up to a multiplicative constant, an antiderivative of the mutation rate.
\bigbreak
{\bf The model}
We assume that a given cancer appears after the following two steps. First we must have a certain mutation and then a clonal expansion of the mutant. Let $N(t)$ be the number mutations up to age $t$. We make the following hypotheses for $N$.

(i) $N(0)=0$.

(ii) The number of mutations in disjoint time intervals are independent random variables. That is, for any integer $k\geq 1$ and any sequence of times $t_0<t_1<\dots<t_k$ we have that the random variables $N(t_n)-N(t_{n-1})$, $1\leq n\leq k$, are independent.

(iii) The probability that there is exactly one mutation during the time interval $(t,t+h)$ is of the order of $\mu(t) h$ as $h$ goes to 0. The mutation rate $\mu\geq 0$ is a function of time. More precisely, we have
$$\lim_{h\to 0}{P(N(t+h)-N(t)=1)\over h}=\mu(t).$$

(iv) The probability that there are two or more mutations during the time interval $(t,t+h)$ is of order less than $h$ as $h$ goes to 0. That is,
$$\lim_{h\to 0}{P(N(t+h)-N(t)\geq 2)\over h}=0.$$

These four hypotheses are quite natural for a mutation process. They model the complete randomness of mutations. Under these hypotheses the process $N$ is necessarily a Poisson process with rate
$$m(t)=\int_0^t\mu(s)ds.$$
See, for instance, Section 5.4.1 in Ross (2006). In particular,
$$P(N(t)=k)=\exp(-m(t)){m(t)^k\over k!}\hbox{ for any integer }k\geq 0.$$
In order to provoke a cancer a given mutation must be successful. This might mean, for instance, replacing the healthy cells of an organ by the mutated cells or a clonal expansion of the mutated cells. Let $\sigma(t)$ be the probability that a mutation appearing at time $t$ will expand. Assume that each mutation expands or not independently of all other mutations and of the process $N$. Under this assumption the number of mutations that will expand, denoted by $C(t)$, is also a Poisson process with rate  
$$M(t)=\int_0^t\sigma(s)\mu(s)ds.$$
For this fact, see (5.9) in Kingman (1993).
Let $I(t)$ be the probability of cancer by time $t$. It seems reasonable to assume that the time it takes for a mutation to expand in the body is small (weeks or months) compared to the time it takes for the appearance of a successful mutation (years or decades). Hence, we may approximate $I(t)$ by $P(C(t)\geq 1)$. 
$$I(t)=P(C(t)\geq 1)=1-\exp(-M(t)).$$
Even at old age the probability of a given cancer is small. Therefore, at all ages $t$, $M(t)$ must be close to 0 and we may approximate
$$I(t)\sim M(t)=\int_0^t\sigma(s)\mu(s)ds.$$
Suppose now that $\sigma$ is a constant (the probability of success for a given mutation does not depend on age). Then the cancer incidence $I$ is up to a multiplicative constant an antiderivative of the mutation rate $\mu$. Hence, in order to model a given incidence $I$ one may define $\mu$ as constant times the derivative of $I$ and then adjust the constant.

For instance, if $I(t)=ct^\gamma$ then set 
$$\mu(s)=\mu_0 s^{\gamma-1}\hbox{ where } \mu_0={c\gamma\over \sigma}.$$

There are several possible interpretations for $\sigma$. It may be the fixation probability (the probability that eventually all cells of the organ will have the mutation) of a mutation in a given organ. Using a Moran model as in Michor et al. (2006) the fixation probability is given by
$$\sigma={1-1/r\over 1-1/r^n}$$
where $r$ is the relative fitness of the mutation and $n$ is the number of cells in the organ.

We get another expression for $\sigma$ if we assume that a mutant expands according to a branching process. Assume that each mutant cell divides with probability $p$ or dies with probability $1-p$ independently of each other and of everything else. Starting with a single mutant the probability of an ever lasting expansion is
$$\sigma=2-{1\over p}\hbox{ for }p>1/2.$$
This expression for $\sigma$ is obtained by observing that
this is a branching process and $1-\sigma$ (the probability that the mutation dies out) is a solution of
$$x=(1-p)+px^2.$$
See, for instance I.9 in Schinazi (1999). 
\bigbreak
{\bf Discussion. }We point out that plausible hypotheses (i)-(iv) the number of mutations in a one-mutation model is a non-homogeneous Poisson process. By neglecting the time of fixation (or expansion) and the time of detection the age incidence for this model is (up to a multiplicative constant) an anti-derivative of the mutation rate. Note that we would get the same result if instead of neglecting these two times we would take them to be constant (as proposed by Moolgavkar and Knudson (1981)). According to our result if the mutation rate is constant the incidence is necessarily linear in age. In contrast, Michor et al. (2006) claim that the incidence may increase as the cube of age for a constant mutation rate in a one-mutation model. They, however, do not neglect the time for a mutant to fixate nor the time to detect the tumor. They take the first time to be constant and the second one to be random.

The idea of using a non-homogeneous Poisson process for carcinogenesis goes back to at least Kendall (1960) although it seems to have been ignored by the recent literature .  
See also Moolgavkar and Venzon (1979), Moolgavkar and Knudson (1981)) and Whittemore and Keller (1978). Our main result, however, seems to be new.

\bigbreak
{\bf References.}

P. Armitage (1985). Multistage models of carcinogenesis. {\it Environmental Health Perspectives}, 63, 195-201.

P. Armitage and R. Doll (1954). The age distribution of cancer and multistage theory of carcinogenesis, {\it British Journal of Cancer}, 8, 1-12.

D.G.Kendall (1960). Birth-and-death processes, and the theory of carcinogenesis. {\it Biometrika}, 47, 13-21.

J.F.C. Kingman (1993). {\it Poisson processes.} Oxford University Press.

F. Michor, Y. Iwasa and M.A. Nowak (2006). The age incidence of chronic myeloid leukemia can be explained by a one mutation model.{\it PNAS}, 103, 14931-14934.

S.H. Mollgavkar and A.G. Knudson Jr. (1981). Mutation and cancer: a model for human carcinogenesis. {\it Journal of the National Cancer institute}, 66, 1037-1051.

S.H. Mollgavkar and D.J. Venzon (1979). Two-event models for carcinogenesis: incidence curves for childhood and adult tumors. {\it Mathematical Biosciences}, 47, 55-77.

S. Ross (2006). {\it Probability Models}, 9th edition, Academic Press.

R.B. Schinazi (1999) {\it Classical and spatial stochastic processes.} Birkhauser.

A. Whittemore and J.B.Keller (1978). Quantitative theories of carcinogenesis. {\it SIAM Review}, 20, 1-30.

\bye